\documentclass[journal,twocolumn]{IEEEtran}
\usepackage{amssymb,amstext,amsmath,amsthm}
\usepackage{multicol}
\usepackage{epsfig,float,lmodern}
\usepackage[compress]{cite}
\usepackage[normalem]{ulem}
\usepackage{flushend}
\usepackage{color}

\usepackage{tcolorbox}
\tcbset{width=0.5\textwidth,boxrule=0pt,colback=red,arc=0pt,auto outer arc,left=0pt,right=0pt,boxsep=5pt}

\newtheorem{definition}{Definition}

\DeclareMathOperator{\sign}{sign}

\begin{document}

\title{
 Solution of the Gardner problem on \\
 the lock-in range of phase-locked loop
}

\author{
  Kuznetsov~N.V., Leonov~G.A., Yuldashev~M.V., Yuldashev~R.V.
\thanks{
}
\thanks{
N.V. Kuznetsov$^{\,a,b,c}$, G.A. Leonov$^{\,b,d}$, M.V. Yuldashev$^{\,b}$,
and R.V. Yuldashev$^{\,b}$ are from
($^a$) Modeling Evolutionary Algorithms Simulation and Artificial Intelligence,
Faculty of Electrical \& Electronics Engineering,
Ton Duc Thang University, Ho Chi Minh, Vietnam;
($^b$) Faculty of Mathematics and Mechanics,
Saint-Petersburg State University, Russia;
($^c$) Dept. of Mathematical Information Technology,
University of Jyv\"{a}skyl\"{a}, Finland;
($^d$) Institute of Problems of Mechanical Engineering RAS, Russia;
(corresponding author email: nikolayv.kuznetsov@tdt.edu.vn, nkuznetsov239@gmail.com).
}
}

\maketitle


\begin{abstract}
The lock-in frequency and lock-in range concepts were introduced in 1966 by Floyd Gardner
to describe the frequency differences of phase-locked loop based circuit
for which the loop can acquire lock within one beat, i.e. without cycle slipping.
These concepts became popular among engineering community and were given
in various engineering publications.
However, rigorous mathematical explanations of these concepts turned out to be a challenging task.
Thus, in the 2nd edition of Gardner's well-known work {\emph{Phaselock Techniques}} he wrote
that {\it{``despite its vague reality, lock-in range is a useful concept}''}
and posed the problem {\it{``to define exactly any unique lock-in frequency}''}.
In this paper an effective solution for Gardner's problem on the definition
of the unique lock-in frequency and lock-in range is discussed.
The lock-in range and lock-in frequency computation
is explained on the example of classical second-order PLL
with lead-lag and active proportional-integral filters.
The obtained results can also be used for the lock-in range computation
of such PLL-based circuits as
two-phase PLL, two-phase Costas loop, BPSK Costas loop,
and optical Costas loop (used in intersatellite communication).
\end{abstract}

\begin{IEEEkeywords}
Gardner problem on lock-in range, cycle slipping, pull-in range,
phase-locked loop, analog PLL, Costas loop,
lead-lag filter, active PI,
cylindrical phase space, separatrix, nonlinear analysis, Lyapunov function.
\end{IEEEkeywords}

\section{Introduction}
\IEEEPARstart{T}{he} phase-locked loop (PLL) is
a \emph{nonlinear control feedback loop} used for synchronization of
the controlled oscillator signal
to the reference oscillator signal (\emph{master-slave synchronization}).
Various PLL-based circuits allow achieving different degrees of synchronization.
For example, the two-phase PLL allows, theoretically, to achieve
\emph{complete synchronization} of signals
(i.e. to get signals with the same frequency and constant phase difference);
classical analog PLL with sinusoidal signal
allows achieving \emph{almost complete synchronization} of signals
(i.e. the difference between phases is almost constant,
the difference between frequencies is almost zero);
PLL with phase-frequency detector is used
to achieve \emph{discrete synchronization}
(i.e. switching points of the signals, e.g. zero crossing points, are synchronized only).
The process of synchronization is called \emph{transient or acquisition process}.
After synchronization is achieved, i.e. transient process is over,
the PLL is said to be in a \emph{locked state}.
The transient process depends on the initial state of the loop
and oscillators' frequencies.
\begin{figure}[!ht]
 \centering
 \includegraphics[width=0.95\linewidth]{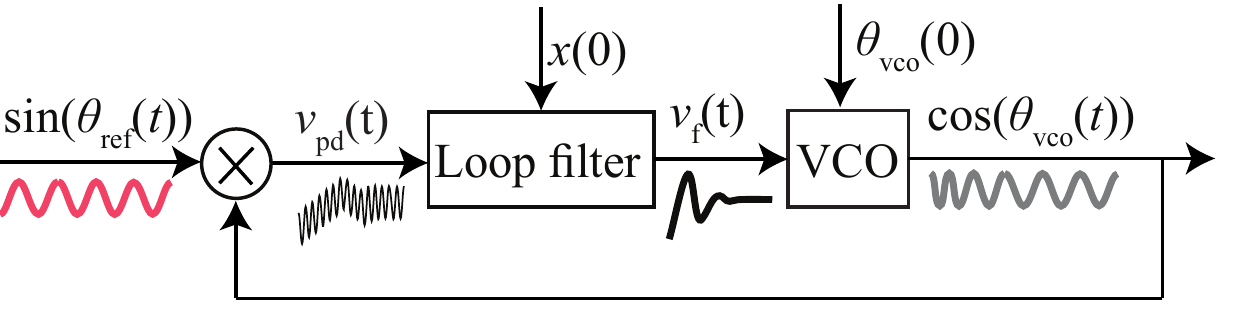}
 \caption{
   Classical PLL with multiplier and sinusoidal signals in the signal space.
 }
 \label{pll-sin-model}
\end{figure}
Consider master-slave synchronization of oscillators by the classical PLL (see, e.g. \cite{Gardner-1966}).
Fig.~\ref{pll-sin-model} shows operation principle of classical PLL, without noise,
in a \emph{signal space} (\emph{the signal space model} of PLL).
Here the input is the signal of the reference oscillator
$\sin(\theta_{\rm{ref}}(t))$ with phase
$\theta_{\rm{ref}}(t)$ and instantaneous frequency
$\dot\theta_{\rm{ref}}(t)=\omega_{\rm{ref}}(t)$.
The signal of the voltage controlled oscillator (VCO)
is $\cos(\theta_{\rm{vco}}(t))$ with phase $\theta_{\rm{vco}}(t)$
and instantaneous frequency
$\dot\theta_{\rm{vco}}(t)=\omega_{\rm{vco}}(t)$.
The multiplier generates a signal
$v_{\rm{pd}}(t)=v_{e}(t)+v_{\rm{hf}}(t)$
which is the sum of the error signal
$v_{e}(t)=\frac{1}{2}\sin(\theta_{\rm{ref}}(t)-\theta_{\rm{vco}}(t))$
and unwanted high-frequency signal
$v_{\rm{hf}}(t)=
\frac{1}{2}\cos(\theta_{\rm{ref}}(t)+\theta_{\rm{vco}}(t))$.
The loop filter suppresses unwanted high-frequency oscillations and
its output $v_{\rm f}(t)$ adjusts the free-running (quiescent) frequency $\omega^{\rm free}_{\rm vco}$ of the VCO
to the frequency of the input signal.
Similarly, operation with square waveform signals
$\sign\sin(\theta_{\rm{ref}}(t))$ and $\sign\cos(\theta_{\rm{vco}}(t))$
can be considered.
Initial state of the loop is $\theta_{\rm vco}(0)$
(initial phase shift of the VCO signal with respect to the reference signal)
and $x(0)$ (initial state of the loop filter).
For the reference signal with constant high frequency
$\omega_{\rm{ref}}(t) \equiv \omega_{\rm{ref}}$ the loop
can achieve only \emph{almost complete synchronization} of the signals
(because of non ideal filtration of high-frequency signal)
when in a locked state
\(
|\omega_{\rm vco}(t) - \omega_{\rm ref}| \leq \omega_{e}^{\rm max}, \
|\theta_e(t) - \theta_{e}(0)| \leq \theta_{e}^{\rm max}, \ t \geq 0,
\)
where constants $\omega_{e}^{\rm max}$ and $\theta_{e}^{\rm max}$
are determined by engineering requirements for a particular application.

Important issues in the design of PLL are
(see, e.g. pioneering monographs \cite{Gardner-1966,ShahgildyanL-1966,Viterbi-1966},
published in 1966, and a rather comprehensive bibliography of pioneering works in \cite{LindseyT-1973}):
estimation of the ranges for deviation between oscillators' frequencies
for which a locked state can be achieved (i.e. the synchronization is possible),
analysis of the locked states stability, and study of possible transient processes.

In \cite{Gardner-1966} Floyd Gardner introduced a lock-in concept:
``\emph{If, for some reason, the frequency difference between input and VCO is less than the loop
bandwidth, the loop will lock up almost instantaneously
without slipping cycles. The maximum frequency difference
for which this fast acquisition is possible is called the lock-in
frequency}''.
The above notion of \emph{lock-in frequency} and
corresponding definition of the \emph{lock-in range}
(called also a \emph{lock range}\cite[p.256]{Yeo-2010-book}, a \emph{seize range} \cite[p.138]{Egan-2007-book})
are given in various engineering publications.\footnote{
See, e.g.
\cite[p.34-35]{Best-1984},\cite[p.161]{Wolaver-1991},\cite[p.612]{HsiehH-1996},\cite[p.532]{Irwin-1997},\cite[p.25]{CraninckxS-1998-book},
\cite[p.49]{KiharaOE-2002},\cite[p.4]{Abramovitch-2002},\cite[p.24]{DeMuerS-2003-book},\cite[p.749]{Dyer-2004-book},\cite[p.56]{Shu-2005},
\cite[p.112]{Goldman-2007-book},\cite[p.61]{Best-2007},\cite[p.138]{Egan-2007-book},
\cite[p.576]{Baker-2011},\cite[p.258]{Kroupa-2012},
\cite[p.387]{Middlestead-2017}
}
However, in general, even for zero frequency difference
there may exist initial states of loop such that
cycle slipping may take place.
Thus, consideration of initial state of the loop is of utmost importance
for the cycle slip analysis and, therefore,
the original concept of lock-in frequency lacks rigor and requires clarification.
In 1979\footnote{A year later, in 1980, F.  Gardner was elected IEEE Fellow
\emph{for contributions to the understanding and applications of phase lock loops}.}
F.~Gardner in the 2nd edition of his well-known work, \emph{Phaselock Techniques},
\cite[p.70]{Gardner-1979-book} (see also the 3rd edition \cite[p.187-188]{Gardner-2005-book} published in 2005) wrote that
 ``{\it{despite its vague reality, lock-in range is a useful concept}}''
 and formulated the following problem:
 ``{\it{there is no natural way to define exactly any unique lock-in frequency}}''.

Recently a rigorous clarification of the lock-in range notion
and solution to Gardner's problem were suggested in
\cite{KuznetsovLYY-2015-IFAC-Ranges,LeonovKYY-2015-TCAS}. 
Below we suggest further extension of the lock-in range notion
for the signal space PLL model.

\begin{definition}[Pull-in and lock-in ranges] 
\label{def-lock-in}
For a certain free-running frequency of the VCO
the largest symmetric interval, around zero\footnote{
In general, for high-order filters
the required behavior can be observed on a union of intervals.
Thus, we are interested in an interval which contains zero,
and it is not necessary defined by the maximum frequency difference
that satisfies the requirements \cite{KuznetsovLYY-2015-IFAC-Ranges,LeonovKYY-2015-TCAS}.
The largest asymmetrical interval containing zero can also be considered.},
of the difference between the reference frequency and the free-running frequency of the VCO,
such that for corresponding frequencies and arbitrary initial state
the loop acquires a locked state,
is called a \emph{pull-in range}.
If, in addition, the loop in a locked state
after any abrupt change of the reference frequency within the interval
acquires a locked state without cycle slipping,
then the interval is called a \emph{lock-in range}.
\end{definition}

Further this definition is explained on the example of classical second-order PLL model
in a signal's phase space with lead-lag and active proportional-integral (PI) filters.
An effective computational method for the lock-in frequency is discussed
and corresponding estimates are obtained.

\begin{figure}[h]
 \centering
 \includegraphics[width=0.95\linewidth]{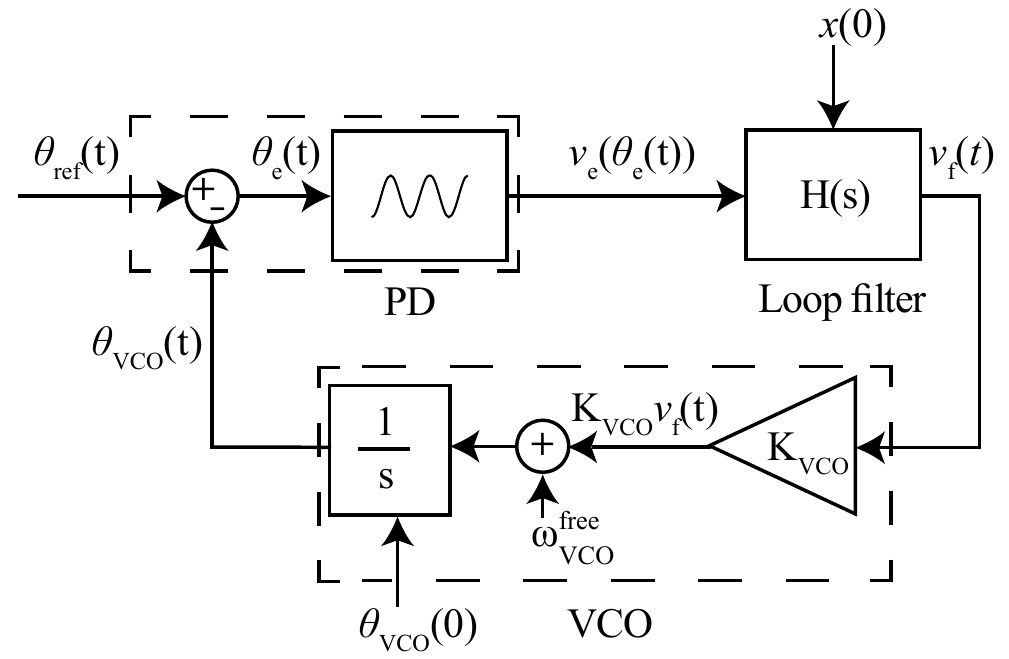}
 \caption{
   Classic analog PLL in the signal's phase space.
 }
 \label{phase-space-model-fig}
\end{figure}

\section{The signal's phase space model} \label{section:model}
Following pioneering engineering publications \cite{Gardner-1966,ShahgildyanL-1966,Viterbi-1966},
consider the operation of classical analog PLL in a \emph{signal's phase space},
where the phase changes of the signals are considered instead of the signals themselves.
The averaging method \cite{KrylovB-1947,KudrewiczW-2007,LeonovKYY-2012-TCASII,LeonovK-2014-book,LeonovKYY-2015-SIGPRO,LeonovKYY-2016-DAN}
allows to reduce the model in the signal space (see Fig.~\ref{pll-sin-model})
to the model in the \emph{signal's phase space} (see Fig.~\ref{phase-space-model-fig})
under certain conditions\footnote{
Remark that their rigorous consideration is often omitted
(see, e.g. classical books \cite[p.12,15-17]{Viterbi-1966},\cite[p.7]{Gardner-1966}),
while their violation may lead to unreliable results
(see, e.g. \cite{PiqueiraM-2003,KuznetsovKLNYY-2015-ISCAS,BestKKLYY-2015-ACC}).
}.
Nowadays nonlinear model in Fig.~\ref{phase-space-model-fig} is widely used
(see, e.g. \cite{Abramovitch-2002,Abramovitch-2004,Best-2007})
to study acquisition processes of various circuits.
In Fig.~\ref{phase-space-model-fig} the phases $\theta_{\rm ref, vco}(t)$ of
the input (reference) and VCO signals
are inputs of the nonlinear phase detector (PD) block. The output of PD is
a function $v_e(\theta_e(t))$ (called PD characteristic),
where
\begin{equation}
  \label{theta_delta_def}
    \theta_e(t) = \theta_{\rm{ref}}(t) - \theta_{\rm{vco}}(t)
\end{equation}
is named the phase error.
For the classical PLL with sinusoidal signals\footnote{
Note that the PD characteristic $v_e(\theta_{e})$ depends
on the waveforms of the considered signals
and can be represented in term of their Fourier coefficients \cite{LeonovKYY-2012-TCASII,LeonovKYY-2015-SIGPRO,KuznetsovLSYY-2015-PD}.
}
and a two-phase PLL\footnote{
A two-phase modification of the classical PLL (see, e.g. \cite{emura2000high,BestKLYY-2014-IJAC,KuznetsovLYY-2017-CNSNS}),
in contrast to the classical one, does not contain high-frequency components:
i.e. $v_{\rm pd}(t) = v_e(\theta_e(t))$.
} we have
\begin{equation}\label{sinusoidal_pd}
  v_e(\theta_e)= \frac{1}{2}\sin(\theta_e),
\end{equation}
in case of square waveforms of reference and VCO signals
$v_{e}(\theta_e)$ has triangular waveform
\begin{equation} \label{triangular_pd}
\begin{aligned}
  & v_e(\theta_e)=\left\{\!\!\!\!
    \begin{array}{ll}
      \frac{2}{\pi}(\theta_e-2\pi), &\!\!\!\!\!\! \rm{if\ } (\theta_e {\, \rm mod\,} 2\pi)
      \in [-\frac{\pi}{2},\frac{\pi}{2}]\\
      2 - \frac{2}{\pi}(\theta_e-2\pi), &\!\!\!\! \rm{ if\ } (\theta_e {\, \rm mod\,} 2\pi)
      \in [\frac{\pi}{2},\frac{3}{2}\pi]
    \end{array}
    \right.
\end{aligned}
\end{equation}

The relationship between the input $v_e(\theta_e(t))$
and the output $v_f(t)$ of the Loop filter is as follows:
\begin{equation}\label{loop-filter}
 \begin{aligned}
 & \dot x = A x + b v_e(\theta_e(t)),
 \ v_{\rm f}(t) = c^*x + hv_e(\theta_e(t)),
 \end{aligned}
\end{equation}
where $A$ is a constant $n\times n$ matrix,
$x(t) \in \mathbb{R}^n$ is the filter state,
$x(0)$ is the initial state of filter,
$b$ and $c$ are constant vectors, and h is a number.
The Loop filter transfer function has the form\footnote{
  In the control theory the transfer function is often defined with the opposite sign
  (see, e.g. \cite{LeonovK-2014-book}): $H(s) = c^*(A-sI)^{-1}b-h.$
}:
\begin{equation}
  H(s) = -c^*(A-sI)^{-1}b+h,
\end{equation}
where $^*$ is matrix transposition, $I$ is $n\times n$ unit matrix.
A lead-lag filter
or a PI filter \cite{Best-2007}
are usually used as the Loop filters.
The control signal $v_f(t)$ adjusts the VCO frequency:
\begin{equation} \label{vco first}
   \dot\theta_{\rm{vco}}(t) = \omega_{\rm{vco}}(t) = \omega_{\rm{vco}}^{\text{free}}
   + K_{\rm vco}v_f(t),
\end{equation}
where $\omega_{\rm{vco}}^{\text{free}}$ is the VCO free-running frequency and $K_{\rm vco}$
is the VCO gain.
Nonlinear VCO models can be similarly considered, see, e.g.
\cite{Margaris-2004,Suarez-2009,BonninCG-2014,BianchiKLYY-2016}.
The frequency of the input signal (reference frequency) is usually assumed
to be constant:
\begin{equation}\label{omega1-const}
  \dot\theta_{\rm{ref}}(t) = \omega_{\rm{ref}}(t) \equiv \omega_{\rm{ref}}.
\end{equation}
The difference between the reference frequency and the VCO free-running frequency
is denoted as $\omega_e^{\text{free}}$:
\begin{equation}
  \label{omega_delta_def}
  \begin{aligned}
    & \omega_e^{\text{free}} \equiv \omega_{\rm{ref}} - \omega_{\rm{vco}}^{\text{free}}.
  \end{aligned}
\end{equation}

By combining equations \eqref{theta_delta_def}, \eqref{loop-filter}, and \eqref{vco first}--\eqref{omega_delta_def}
a \emph{nonlinear mathematical model in the signal's phase space} is obtained
(i.e. in the state space: the filter's state $x$ and the difference between the signal's phases $\theta_e$):
\begin{equation}\label{final_system}
\begin{aligned}
   & \dot{x} = A x + b v_e(\theta_{e}), \\
   & \dot\theta_{e} = \omega_{e}^{\text{free}}
   - K_{\rm vco} \big(c^*x + hv_e(\theta_{e})\big).
 \end{aligned}
\end{equation}

Classical PD characteristics are bounded piecewise-smooth $2\pi$ periodic\footnote{
If $v_e(\theta_e)$ has another period
(e.g. $\pi$ for the Costas loop models),
it has to be considered in the further discussion instead of $2\pi$.
} functions, thus it is convenient to assume that ($\theta_e$ mod $2\pi$) is a cyclic variable,
and the analysis is restricted to the range of $\theta_e(0) \in [-\pi, \pi)$.

System  \eqref{final_system} with an odd PD characteristic (i.e. $v_{e}(-\theta_e)=-v_{e}(\theta_e)$) is not changed by the transformation
\begin{equation}\label{odd-change}
  \big(\omega_{e}^{\text{free}},x(t),\theta_{e}(t)) \rightarrow
  \big(-\omega_{e}^{\text{free}},-x(t),-\theta_{e}(t)),
\end{equation}
and it allows to study system \eqref{final_system} for $\omega_e^{\text{free}}>0$ only,
introducing the concept of \emph{frequency deviation}:
\begin{equation}\label{eq:fd}
  |\omega_e^{\text{free}}| = |\omega_{\rm{ref}} - \omega_{\rm{vco}}^{\text{free}}|.
\end{equation}

\subsection{Pull-in and lock-in ranges}
For nonlinear mathematical model in the signal's phase space \eqref{final_system}
we can consider the conditions of complete synchronization\footnote{
If necessary conditions for the averaging are satisfied
(i.e. the considered frequencies are sufficiently large)
then complete synchronization
for the mathematical model in the signal's phase space implies
almost complete synchronization for the mathematical model in the signal space
\cite[p.88]{mitropolsky1967averaging},\cite{KudrewiczW-2007}.
},
i.e. the frequency error is zero and the phase error is constant:
\begin{equation} \label{steady-phase}
 \dot\theta_e(t)\equiv 0, \theta_e(t) \equiv \theta_s.
\end{equation}
For most of the considered loop filters  (i.e. controllable and observable)
the above equation implies that filter state is also constant:
\begin{equation}\label{steady-filter}
   x(t)\equiv x_s.
\end{equation}
Thus, in complete synchronization
the locked states of the model in the \emph{signal's phase space}
are the equilibrium points of system \eqref{final_system}.
The equilibrium points can be considered as
a multiple-valued function of variable $\omega_e^{\text{free}}$:
$\big(x_s(\omega_e^{\text{free}}), \theta_s(\omega_e^{\text{free}})\big)$.
Note, that for all practically used phase detectors the characteristics $v_e(\theta_e)$ (e.g. sinusoidal and triangular) have exactly one stable and one unstable equilibriua on each period.



An important characteristic of PLL is the set of $\omega_{e}^{\rm free}$
such that the model acquires locked state for any initial state.
\begin{definition}[Pull-in range of the \emph{signal's phase space} model,
see \cite{KuznetsovLYY-2015-IFAC-Ranges,LeonovKYY-2015-TCAS,BestKLYY-2016}]
The largest interval of frequency deviations
$|\omega_{e}^{\rm free}| \in [0,\omega_{\text{pull-in}})$
such that the \emph{signal's phase space} model \eqref{final_system} acquires a locked state
for arbitrary initial state $(x(0),\theta_e(0))$ is called
a  \emph{pull-in range}, $\omega_{\text{pull-in}}$ is called a \emph{pull-in frequency}.
\end{definition}


\begin{definition}[Cycle slipping]
Let PD characteristic $v_e(\theta_e)$ be $2\pi$-periodic function.
If
\begin{equation}\label{eq-cs-sup}
   \sup\limits_{t>0} |\theta_{e}(0) - \theta_{e}(t)| \geq 2\pi
\end{equation}
it is then said that \emph{cycle slipping} occurs
\cite[p.131]{Stensby-1997},\cite{LeonovKYY-2015-TCAS}.
Sometimes $\limsup\limits_{t\to+\infty}$
is considered instead of $\sup\limits_{t>0}$ in \eqref{eq-cs-sup}.
\end{definition}
Further we consider only cycle slipping caused by deterministic dynamics \cite{Stoker-1950,Gardner-1966,LeonovRS-1992,SmirnovaP-2016-ECC},
noise-induced cycle slips are studied, 
e.g. in \cite{Tikhonov-1959,Stensby-1997,Talbot-2012-book}.

\begin{figure}[h]
 \centering
 \includegraphics[width=0.95\linewidth]{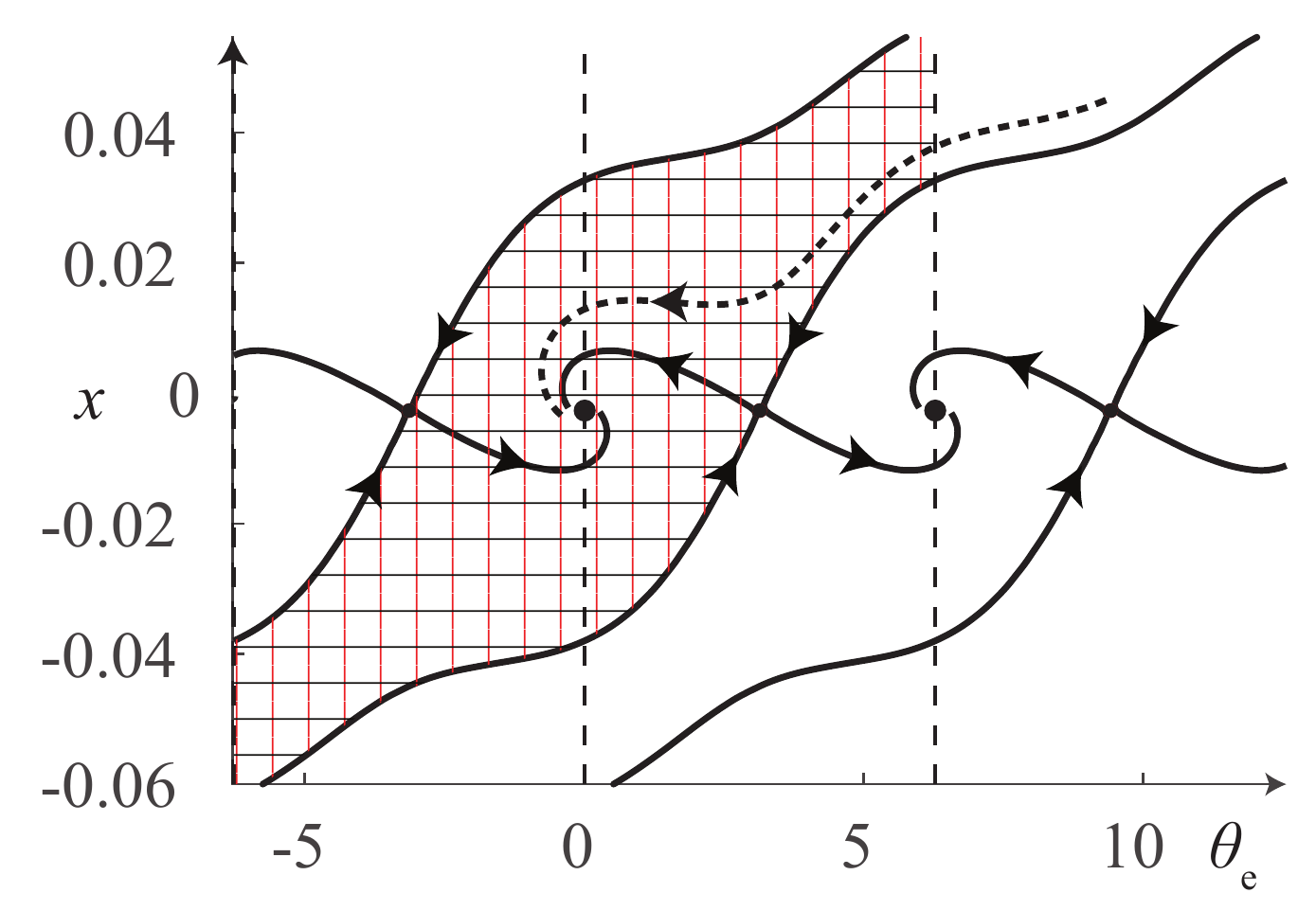}
 \caption{
 Phase portrait of the signal's phase space model of
 the classical PLL with lead-lag filter:
 $H(s)= \frac{1+s\tau_2}{1+s\tau_1}$,
 $\tau_1 = 6.33\cdot10^{-2}$,
 $\tau_2 = 1.85\cdot10^{-2}$,
 $K_{\rm vco}=250$, $v_e(\theta_e)= \frac{1}{2}\sin(\theta_e)$, $\omega_e^{\text{free}} = 0$.
 Dots are equilibria. Separatrices pass in and out of unstable saddle equilibria,
 outgoing separatrices tend to stable equilibria.
 The local \emph{lock-in domain}, where trajectories tend
 to zero equilibrium without slipping cycles, is shaded.
 The dashed trajectory starts outside the \emph{lock-in domain} and skips cycles.
 }
 \label{lock-in-0}
\end{figure}

The pull-in (or acquisition) process
may take more than one beat note, i.e
the VCO frequency will be slowly tuned
toward the reference frequency and cycle slipping may take place.
Remark that even for zero frequency deviation ($\omega_e^{\text{free}}=0$)
and a sufficiently large initial state of filter $x(0)$
the cycle slipping may take place (see, e.g. dashed trajectory in Fig.~\ref{lock-in-0}).
Thus, consideration of all state variables is of utmost importance
for the cycle slipping analysis.
A numerical study of cycle slipping for the  classical PLL can be found in \cite{AscheidM-1982}.
Analytical tools for estimating the number of cycle slips
 can be found, e.g. in \cite{ErshovaL-1983,LeonovRS-1992,LeonovK-2014-book}.

\begin{definition}[Lock-in range of the \emph{signal's phase space} model,
see \cite{KuznetsovLYY-2015-IFAC-Ranges,LeonovKYY-2015-TCAS,BestKLYY-2016}]
\label{def:lockin-sps}
The largest interval of frequency deviations from the pull-in range:
$|\omega_{e}^{\rm free}| \in [0,\omega_{\text{lock-in}}) \subset [0,\omega_{\text{pull-in}})$,
is called a \emph{lock-in range}
if the signal's phase space model \eqref{final_system}, being in a locked state,
after any abrupt change of $\omega_{e}^{\rm free}$
within the interval acquires a locked state without cycle slipping.
\end{definition}
Remark that in contrast to the pull-in and lock-in ranges
for the signal space model (see Definition~1, Fig.~\ref{pll-sin-model})
the ranges for the signal's phase space model (Fig.~\ref{phase-space-model-fig})
depend only on the difference between frequencies
$|\omega_{\rm{ref}} - \omega_{\rm{vco}}^{\text{free}}|$,
but not on the frequencies $\omega_{\rm{ref}}$ and $\omega_{\rm{vco}}^{\text{free}}$ themselves.
However, it must be remembered that
reduction of the signal space model to
the signal's phase space model by the averaging
is reliable only for sufficiently high frequencies.

Consider the behavior of trajectories (transient processes)
of system \eqref{final_system} in the state space $(x, \theta)$
(see, e.g., the state space for a second-order model in Fig.~\ref{lock-in-0}).
Each of the equilibrium states (locked states)
has a vicinity in the state space,
called a \emph{local lock-in domain},
where corresponding trajectories tend to the locked state without slipping cycles
(see, e.g. shaded local lock-in domain in Fig.~\ref{lock-in-0}
defined by corresponding ingoing separatrices of unstable equilibria).
\emph{The lock-in domain}\footnote{
In \cite[p.50]{Viterbi-1966}) the lock-in domain is called a \emph{frequency lock},
in \cite[p.132]{Stensby-1997},\cite[p.355]{Meyer-2004-book})
the \emph{lock-in range} notion is used to denote the \emph{lock-in domain}.
} is the union of local lock-in domains,
each of which corresponds to one of the equilibria.
The shape of the lock-in domain may vary significantly depending on $\omega_{e}^{\text{free}}$.
We denote by ${\rm D}_{\text{lock-in}}(\omega_{e}^{\text{free}})$ lock-in domain corresponding
to a specific frequency difference $\omega_{e}^{\text{free}}$.
The above Definition~\ref{def:lockin-sps} requires that
the intersection of lock-in domains
 \(
  {\rm D}_{\text{lock-in}}\big((-\omega_\text{lock-in},\omega_\text{lock-in})\big) =
  \bigcap\limits_{|\omega_e^{\text{free}}| < \omega_\text{lock-in}} {\rm D}_{\text{lock-in}}(\omega_{e}^{\text{free}})
 \)
 contains all corresponding equilibria
 \(
\big( x_{s}(\omega_{e}^{\text{free}}),
     \theta_{s}(\omega_{e}^{\text{free}})\big).
\)
Thus, after an abrupt change of frequencies,
the model acquires a locked state without slipping cycles.

\section{Lock-in range of the classical analog PLL}

To estimate numerically the lock-in range for
the signal's phase space model in Fig.~\ref{phase-space-model-fig},
we can use, e.g., direct numerical integration
of ordinary differential equation \eqref{final_system} in MATLAB\footnote{See e.g.
MATLAB ODE solvers ode45 and ode15s.
}
or simulation of the signal's phase space model
in MATLAB Simulink\footnote{
The MATLAB Simulink block \emph{Transfer Fcn (with initial states)}
allows taking into account the initial filter state $x(0)$;
the initial phase error $\theta_\Delta(0)$ can be taken into account
by the property \emph{initial data} of the \emph{Intergator} blocks.
Corresponding initial states in SPICE
(e.g. capacitor's initial charge) 
are zero by default
but can be changed manually \cite{BianchiKLYY-2015,KuznetsovLYY-2017-CNSNS}.
}.

With the usual numerical simulation or integration of dynamical model
we can only observe transient processes (trajectories)
toward the asymptotically stable locked states (equilibria),
i.e. which are maintained under small perturbations
of the phase difference and filter state;
transient processes toward unstable locked states
(see, e.g. ingoing separatrices of saddle equilibria in Fig.~\ref{lock-in-0})
cannot be obtained due to computational errors
(caused by finite precision arithmetic
and numerical approximation of ordinary differential equation solutions)
and do not occur in real systems because of noise.

The usual numerical simulation
allows observing transient processes (trajectories)
toward the asymptotically stable locked states (equilibria) only,
i.e. locked states which are maintained under small perturbations
of the phase difference and filter state.
Transient processes toward  the unstable locked states
(see, e.g. ingoing separatrices of saddle equilibria in Fig.~\ref{lock-in-0})
cannot be observed due to computational errors
(caused by finite precision arithmetic
and numerical approximation of ordinary differential equation solutions)
and do not occur in real systems because of noise.

\begin{figure}[!htp]
  \centering
  \includegraphics[width=0.95\linewidth]{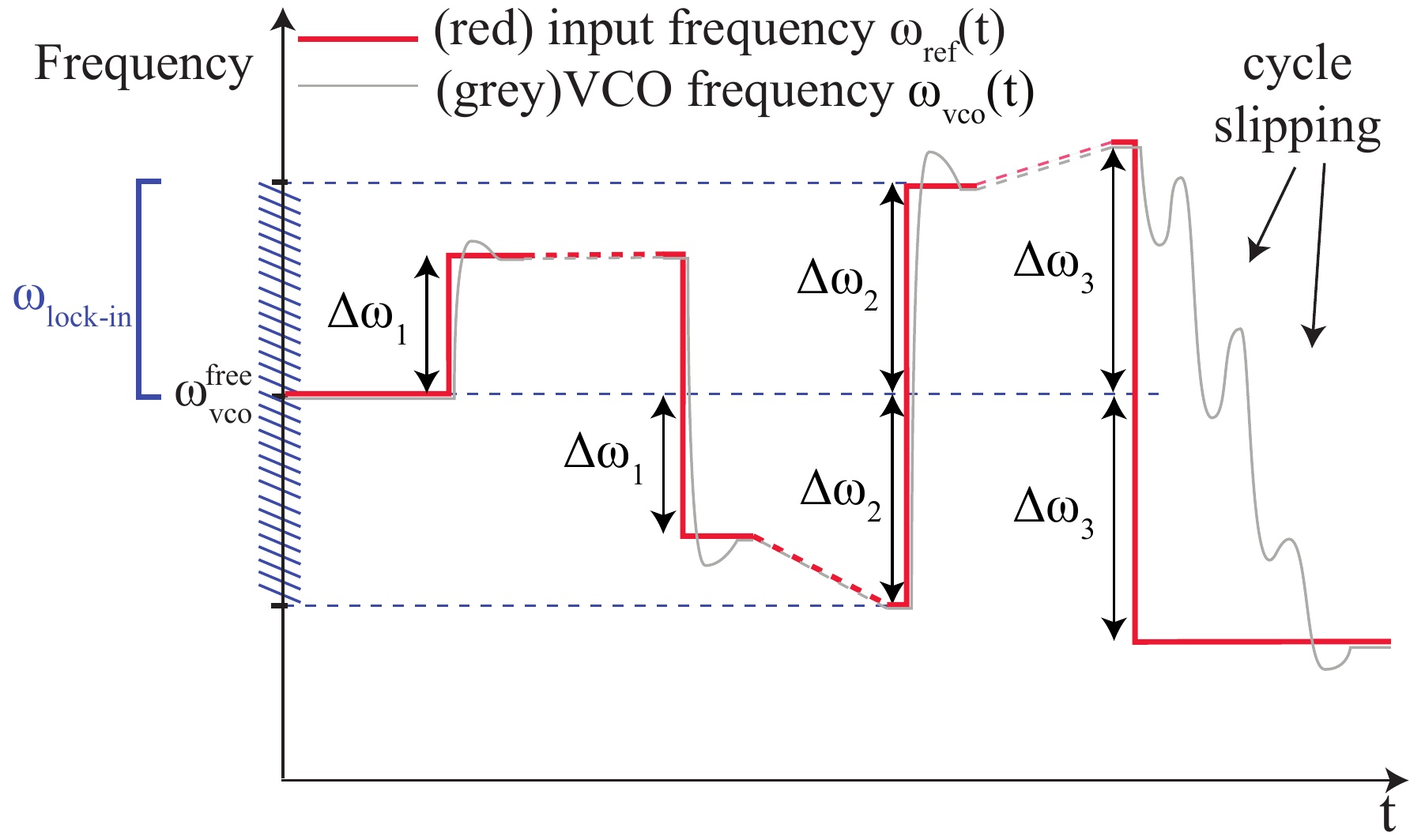}
  \caption{
  Estimation of the lock-in range.}
\label{lock-in-computation-algo}
\end{figure}

Without loss of generality we can fix $\omega_{\rm vco}^{\rm free}$ and vary $\omega_{\rm ref}$.
Assume that the zero frequency deviation $\omega_{e}^{\text{free}}=0$
(i.e. the input frequency $\omega_{\rm{ref}}=\omega_{\rm vco}^{\rm free}$)
belongs to the pull-in range and
the model acquires an asymptotically stable locked state $(\theta_{s}^0,x_{s}^0)$
(without loss of generality we can consider $\theta^0_{s}=0$).
First we abruptly increase the input frequency by sufficiently small
frequency step $\Delta\omega>0$ (i.e. the input frequency becomes
$\omega_{\rm ref} = \omega_{\rm vco}^{\rm free}+\Delta\omega$)
and observe whether corresponding transient process converges to
a locked state without cycle slipping (see Fig.~\ref{lock-in-computation-algo}).
If $|\Delta\omega|$ belongs to the pull-in range then
the model acquires new asymptotically stable locked state
$(\theta^{+}_{s},x^{+}_{s})$
(otherwise we reduce $\Delta\omega$).
If initial locked state belongs to the corresponding lock-in domain:
$(\theta^0_s,x^0_s) \in {\rm D}_{\rm lock-in}(\Delta\omega)$,
then the transient process from $(\theta^0_s,x^0_s)$
converges to the locked state $(\theta^{+}_{s},x^{+}_{s})$ without cycle slipping.
Then we abruptly decrease the input frequency by $2\Delta\omega$, (i.e.
the input frequency becomes $\omega_{\rm ref} = \omega_{\rm vco}^{\rm free}-\Delta\omega$).
Since $|\Delta\omega|$ belongs to the pull-in range
the model acquires new asymptotically stable locked state $(\theta^{-}_{s},x^{-}_{s})$.
If the transient process converges from $(\theta^{+}_s,x^{+}_s)$
to the locked state $(\theta^{-}_{s},x^{-}_{s})$
without cycle slipping,
(i.e. both symmetric locked states belong to the intersection of the lock-in domains:
$(\theta_{s}^{\pm}, x_{s}^{\pm})
\in
{\rm D}_{\rm  lock-in}(|\Delta\omega|) =
{\rm D}_{\rm  lock-in}(+\Delta\omega)\bigcap{\rm D}_{\rm  lock-in}(-\Delta\omega)$),
then we estimate the lock-in range as $[0,\Delta\omega]$.
\begin{figure*}[h]
  \centering
  \includegraphics[width=0.95\linewidth]{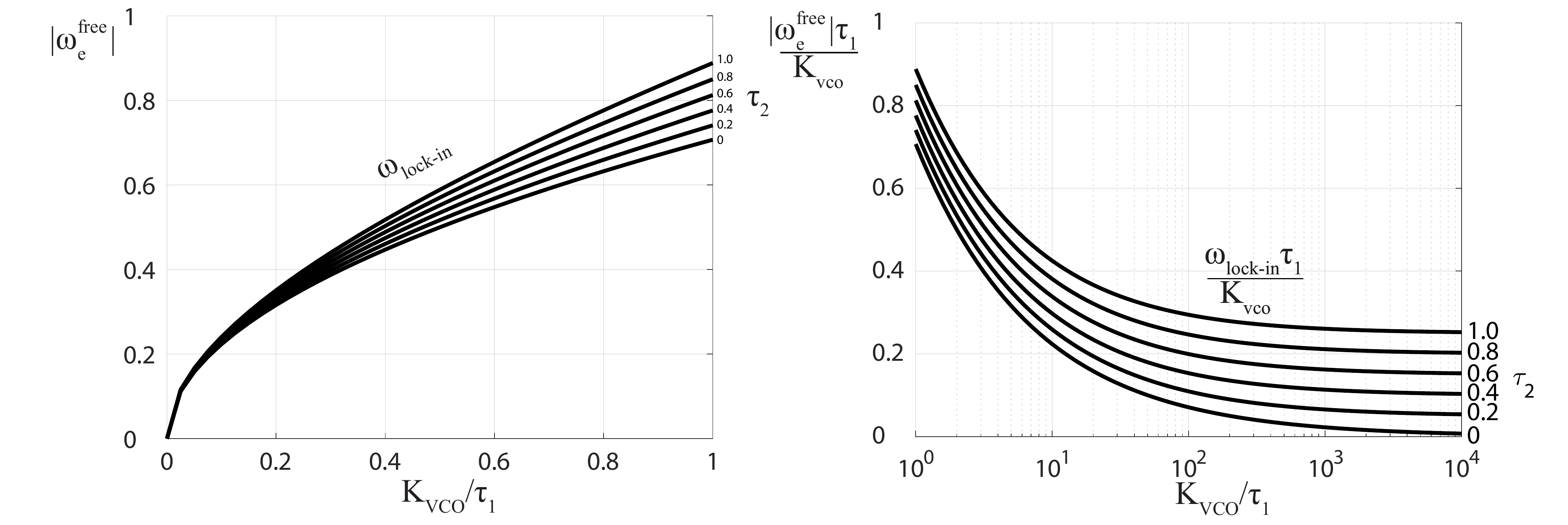}
  \caption{
  Lock-in frequency diagrams for model \eqref{PLLPI} (active PI filter)
  with sinusoidal phase detector characteristics \eqref{sinusoidal_pd}.
  Left subfigure: $\omega_{\text{lock-in}}$ for $0 < K_{\rm vco}/\tau_1 < 1$ in linear scale
  ($\omega_{\text{lock-in}}$ tends to $0$ as $K_{\rm vco}/\tau_1$ tends to 0).
  Right subfigure: $\frac{\omega_{\text{lock-in}}\tau_1}{K_{\rm vco}}$
  for $K_{\rm vco}/\tau_1 > 1$ in logarithmic scale.
  }
  \label{lock-in-perfect-pi-sin}
\end{figure*}

\begin{figure*}[h]
  \centering
  \includegraphics[width=0.95\linewidth]{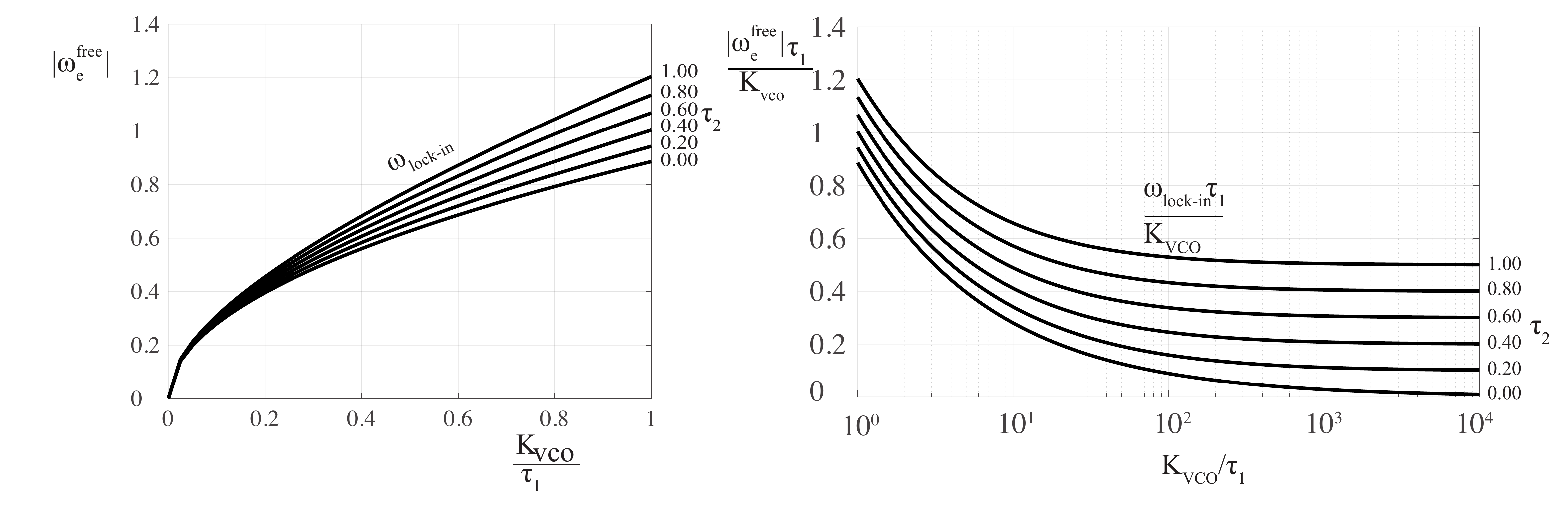}
  \caption{
  Lock-in frequency diagrams for model \eqref{PLLPI} (active PI filter) with
  triangular phase detector characteristics \eqref{triangular_pd}.
  Left subfigure: $\omega_{\text{lock-in}}$ for $0 < K_{\rm vco}/\tau_1 < 1$ in linear scale
  ($\omega_{\text{lock-in}}$ tends to $0$ as $K_{\rm vco}/\tau_1$ tends to 0).
  Right subfigure: $\frac{\omega_{\text{lock-in}}\tau_1}{K_{\rm vco}}$
  for $K_{\rm vco}/\tau_1 > 1$ in logarithmic scale.
  }
  \label{lock-in-perfect-pi-tri}
\end{figure*}
Here we suppose that for a specific $\omega_{\rm ref}$
only one asymptotically stable locked state exist within a periodic interval:
$(\theta_{s}^{\rm st}(\omega_{\rm ref}),x_{s}^{\rm st}(\omega_{\rm ref}))$,
and for all $|\omega_{\rm ref}|<|\Delta\omega|$ we have
  $(\theta_{s}^{\rm st}(\omega_{\rm ref}),x_{s}^{\rm st}(\omega_{\rm ref}))
  \in {\rm D}_{\rm  lock-in}(|\Delta\omega|).$
This allows investigating cycle slipping
for a consecutive increases of frequency deviation $|\omega_{\rm ref}|$
only. Thus, simulation of all possible abrupt changes of
frequency deviation within the lock-in range
(see Fig.~\ref{lock-in-computation-algo}) is unnecessary.


The algorithm for computation of lock-in ranges used above doesn't guarantee
that the model acquires locked state for any initial filter state,
i.e. lock-in range is not inside pull-in range automatically.
To verify the global stability and estimate the pull-in range one can
use the phase plane analysis
(see, e.g. \cite{AndronovVKh-1937,Gubar-1961,Shakhtarin-1969})
or a special modification of direct Lyapunov method
for cylindrical phase space
(see, e.g. \cite{GeligLY-1978,LeonovRS-1992,LeonovK-2014-book,LeonovKYY-2015-TCAS}).
Some examples of the lock-in range computation
by the phase plane analysis
are presented in \cite{LeonovKYY-2015-TCAS,BlagovKKLYY-2016-IFAC,AleksandrovKLNYY-2016-IFAC,AleksandrovKLYY-2016-arXiv-sin,AleksandrovKLYY-2016-arXiv-impulse}.

Below we provide complete solution
of the Gardner problem
on lock-in range computation for the classical
second-order PLL
with lead-lag and active proportional-integral (PI) filters,
and sinusoidal and triangular PD characteristics.

\subsection{Active PI filter}
\begin{figure*}[h]
  \centering
  \includegraphics[width=0.95\linewidth]{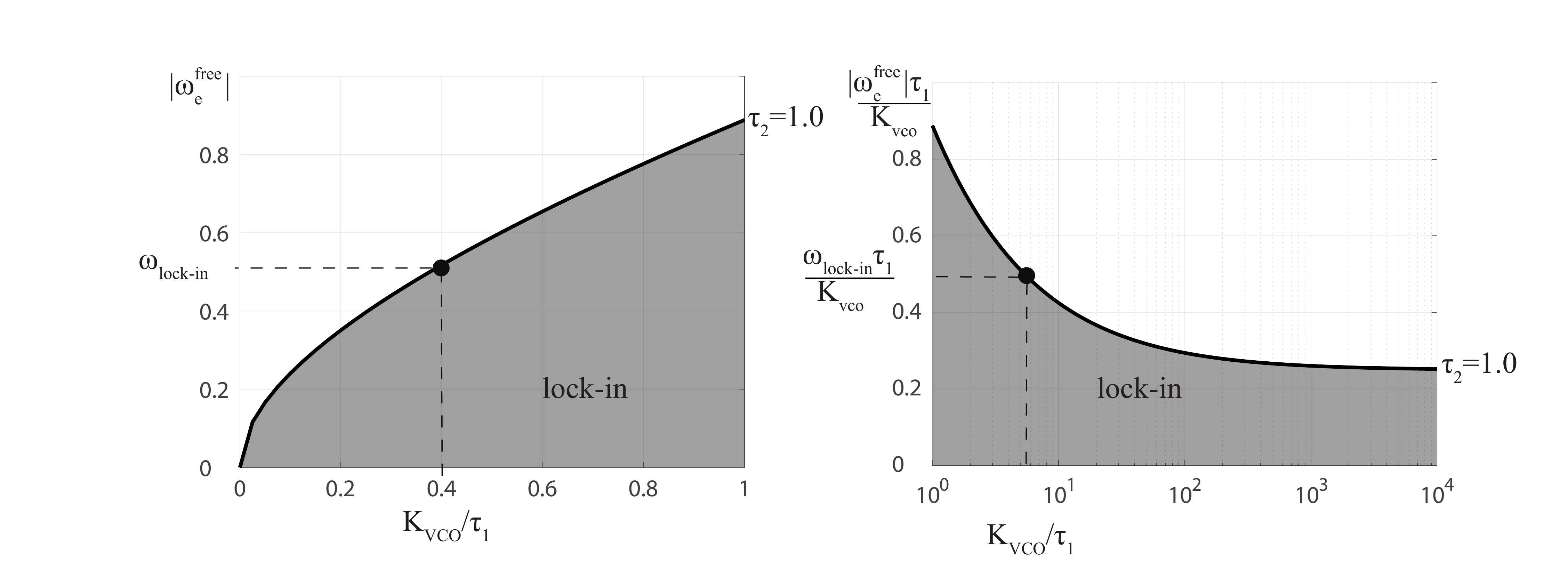}
  \caption{
  Lock-in range depending on $\frac{K_{\rm vco}}{\tau_1}$ (shaded domains)
  for model \eqref{PLLPI} (active PI filter)
  with sinusoidal phase detector characteristics \eqref{sinusoidal_pd} and $\tau_2=1$.
  Left subfigure: $\omega_{\text{lock-in}}$ for $0 < K_{\rm vco}/\tau_1 < 1$ in linear scale.
  Right subfigure: $\frac{\omega_{\text{lock-in}}\tau_1}{K_{\rm vco}}$
  for $K_{\rm vco}/\tau_1 > 1$ in logarithmic scale.
  }
  \label{lock-in-diagram-expl}
\end{figure*}
\begin{figure*}[h]
  \centering
  \includegraphics[width=0.95\linewidth]{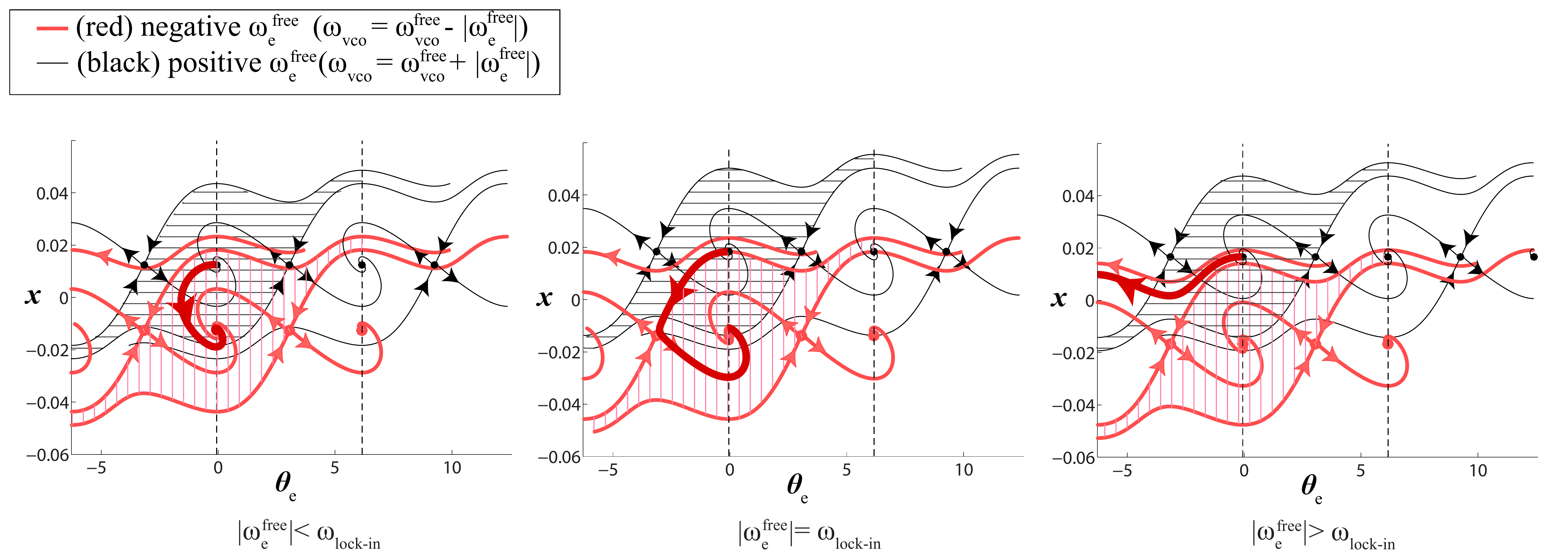}
  \caption{
  Phase portraits for model \eqref{PLLPI} with the following parameters:
  $H(s)= \frac{1+s\tau_2}{s\tau_1}$,
 $\tau_1 = 0.0633$,
 $\tau_2 = 0.0225$,
 $K_{\rm vco}=250$, $v_e(\theta_e)= \frac{1}{2}\sin(\theta_e)$.
 Black color is for the model with positive $\omega_e^{\text{free}}=|\widetilde{\omega}|$.
 Red is for the model with negative $\omega_e^{\text{free}}=-|\widetilde{\omega}|$.
 Equilibria (dots), separatrices pass in and out of the saddles equilibria,
 local lock-in domains are shaded
 (upper black horizontal lines are for $\omega_e^{\text{free}}>0$,
 lower red vertical lines are for $\omega_e^{\text{free}}<0$).
 Left subfig: $\omega_e^{\text{free}} = \pm 50$;
 middle subfig: $\omega_e^{\text{free}} = \pm 60$;
 right subfig: $\omega_e^{\text{free}} = \pm 66$.
 }
 \label{lock-in-pi-pp-fig}
\end{figure*}

 Consider the PLL with the active PI loop filter (called ``perfect PI'')
 having transfer function
 $H(s) = \frac{1+s\tau_2}{s\tau_1}
 $, where $\tau_1,\tau_2 >0$. The model \eqref{final_system} becomes
\begin{equation}\label{PLLPI}
  \begin{aligned}
    & \dot x = \frac{v_e(\theta_e)}{\tau_1},\\
    & \dot\theta_{e} = \omega_{e}^{\rm{free}}
   - K_{\rm{vco}}
   \bigg(
    x + \frac{\tau_2}{\tau_1}v_e(\theta_e)
   \bigg).
  \end{aligned}
\end{equation}
 For PD characteristics \eqref{sinusoidal_pd} and \eqref{triangular_pd}
 the model has one asymptotically stable $(0,-\frac{\omega_e^{\rm{free}}}{K_{\rm vco}})$
 and one unstable $(\pi,-\frac{\omega_e^{\rm{free}}}{K_{\rm vco}})$
 equilibrium points, which repeat on intervals of length $2\pi$.
 Note, that by the change of variables $(x \to x+\frac{\omega_e^{\rm{free}}}{K_{\rm vco}})$
 in system \eqref{PLLPI} we get $\omega_{e}^{\rm{free}} \to 0$.
 Thus, a change of the parameter $\omega_{e}^{\text{free}}$ shifts the phase portrait vertically
 (in the variable $x$) without distorting trajectories,
 which simplifies the analysis of the lock-in domain and range
(see Fig.~\ref{lock-in-pi-pp-fig}).
 The pull-in range of model \eqref{PLLPI} is infinite\footnote{
The rigorous analytical proof can be effectively achieved
by the \emph{special modification of direct Lyapunov method
for cylindrical phase space}
and considering Lyapunov function \cite{Bakaev-1963,LeonovK-2014-book,AlexandrovKLNS-2015-IFAC,LeonovKYY-2015-TCAS}:
$V(x, \theta_e) =
\frac{K_{\rm vco}\tau_1}{2} \left(x - \frac{\omega_e^{\rm free}}{K_{\rm vco}}\right)^2
+ \int \limits_{0}^{\theta_e} v_{e}(s) ds \geq 0$
and $\dot V(x,\theta_e) = -K_{\rm vco}\frac{\tau_2}{\tau_1}v_e^2(\theta_e) \leq 0$.
Here it is important that for any $\omega_e^{\text{free}}$
the set $\dot V(x,\theta_e)\equiv0$ does not contain the whole trajectories of the system except for equilibria.
Note that considered Lyapunov function $V(x,\theta_e)$ is periodic in $\theta_e$
and is bounded for any $||(0,\theta_e)|| \to +\infty$.
Thus, the classical Krasovskii--LaSalle principle on global stability,
where the function has to be radially unbounded
(i.e. $V(x,\theta_e) \to +\infty$ as $||(x,\theta_e)|| \to +\infty$),
cannot be applied directly.}.

 For an arbitrary $\tau_1>0$ the following change of variables
 \(
    x \to x/\tau_1
 \)
 reduces system \eqref{PLLPI} to the form
\begin{equation}\label{PLLPI-reduced-params}
  \begin{aligned}
    & \dot x = v_e(\theta_e),\\
    & \dot \theta_{e} = \omega_{e}^{\rm{free}}
   -  \frac{K_{\rm{vco}}}{\tau_1}
   \bigg(
    x + \tau_2v_e(\theta_e)
   \bigg).
  \end{aligned}
\end{equation}

\begin{figure*}[h]
  \centering
  \includegraphics[width=0.95\linewidth]{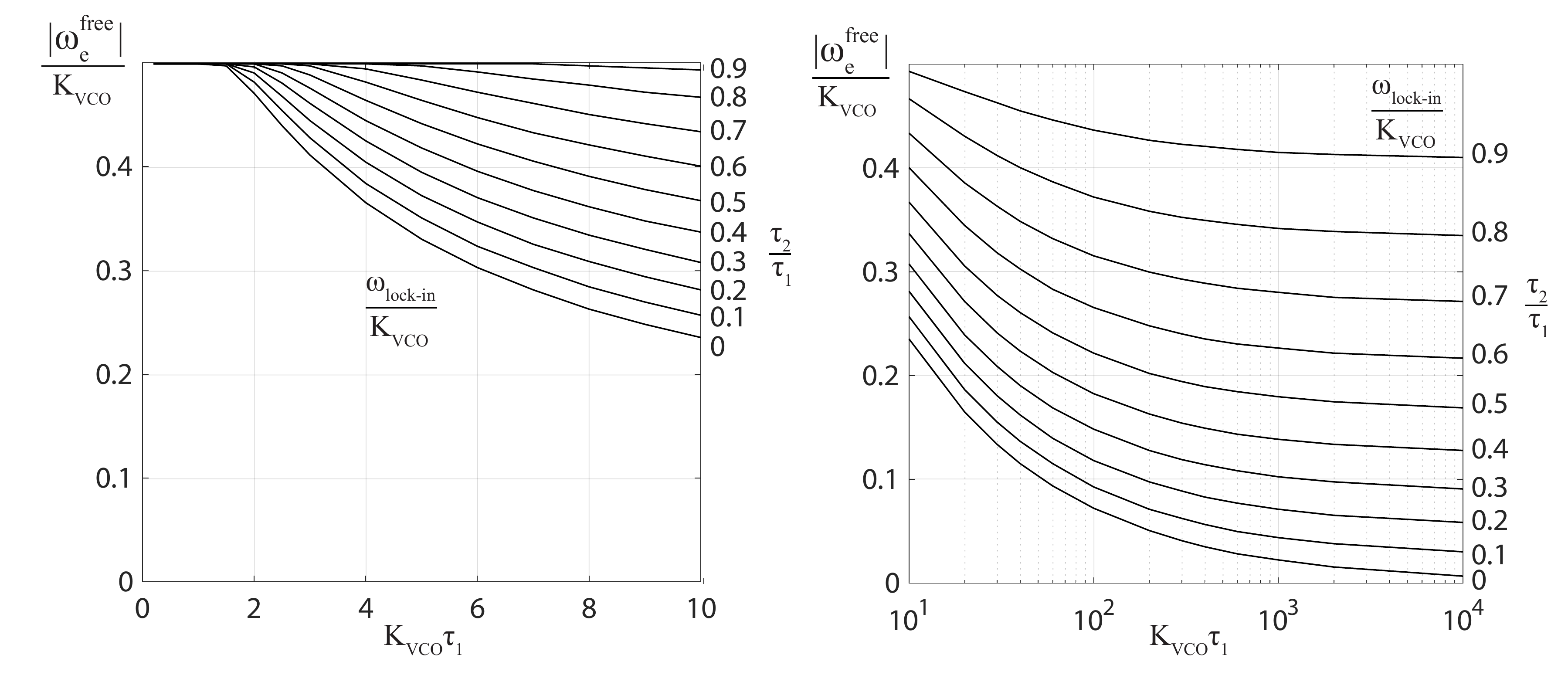}
  \caption{
  Lock-in frequency diagrams for model \eqref{PLLleadlag} (lead-lag filter)
  with sinusoidal phase detector characteristics \eqref{sinusoidal_pd}.
  Left subfigure: $\frac{\omega_{\text{lock-in}}}{K_{\rm vco}}$ for $0 < K_{\rm vco}\tau_1 < 10$
  in linear scale
  ($\frac{\omega_{\text{lock-in}}}{K_{\rm vco}}$ tends to $1/2$ as $K_{\rm vco}\tau_1$ tends to 0).
  Right subfigure: $\frac{\omega_{\text{lock-in}}}{K_{\rm vco}}$ for $K_{\rm vco}\tau_1 > 10$ in logarithmic scale.
  }
 \label{lock-in-diagram-lead-lag-sin}
 \end{figure*}
 \begin{figure*}[h]
  \centering
  \includegraphics[width=0.95\linewidth]{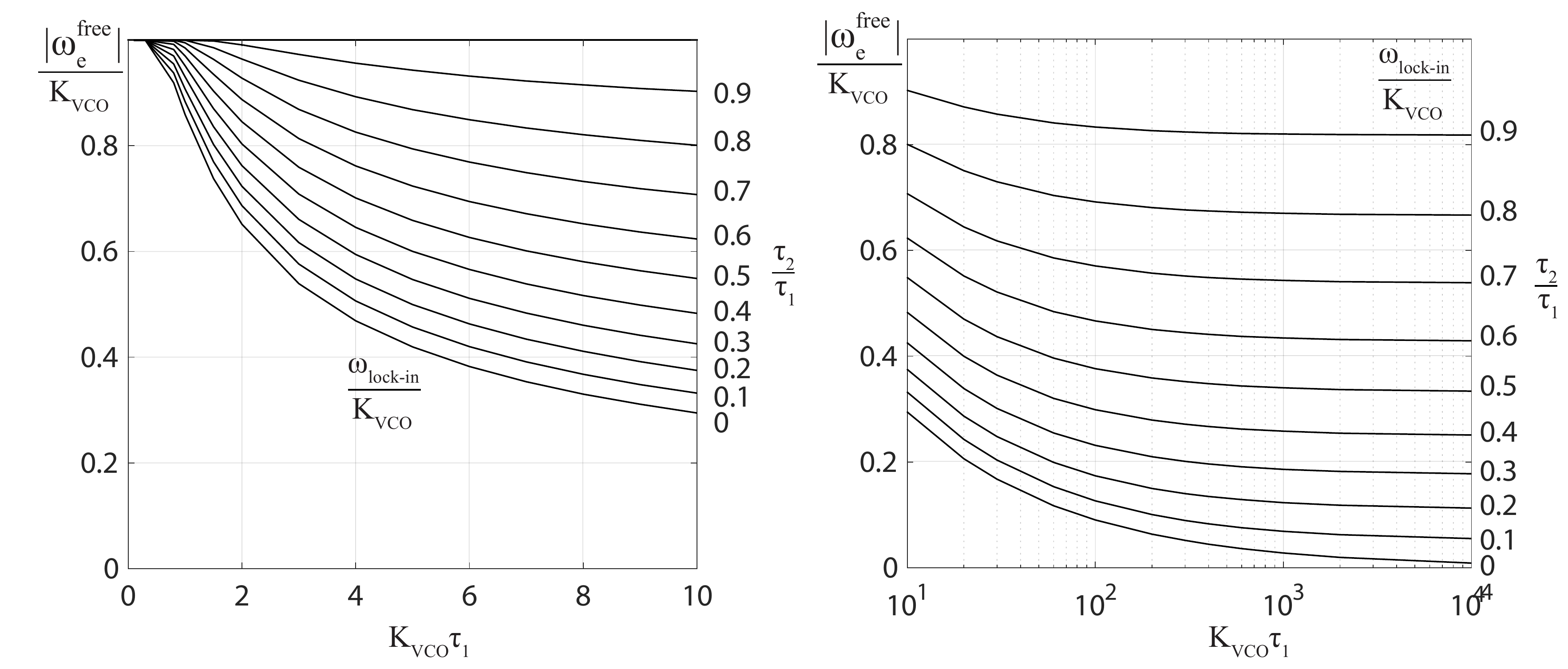}
  \caption{
  Lock-in frequency diagrams for model \eqref{PLLleadlag} (lead-lag filter)
  with triangular phase detector characteristics \eqref{triangular_pd}.
  Left subfigure: $\frac{\omega_{\text{lock-in}}}{K_{\rm vco}}$ for $0 < K_{\rm vco}\tau_1 < 10$
  in linear scale
  ($\frac{\omega_{\text{lock-in}}}{K_{\rm vco}}$ tends to $1$ as $K_{\rm vco}\tau_1$ tends to 0).
  Right subfigure: $\frac{\omega_{\text{lock-in}}}{K_{\rm vco}}$ for $K_{\rm vco}\tau_1 > 10$ in  logarithmic scale.
  }
  \label{lock-in-diagram-lead-lag-tri}
\end{figure*}

This change of variables
rescales the phase portrait of the system (see Fig.~\ref{lock-in-0})
and, thus, does not affect the stability and cycle slipping properties
of the trajectories.
Therefore, $\omega_{\text{lock-in}}$ can be found as a function of
$\tau_2,K_{\rm vco}/\tau_1$.
For model \eqref{PLLPI}
with phase detector characteristics \eqref{sinusoidal_pd} and \eqref{triangular_pd}
the lock-in frequency diagrams
as function of $\tau_2$ (for $\tau_2 \in [0,1]$) and $K_{\rm vco}/\tau_1$
are shown in Fig.~\ref{lock-in-perfect-pi-sin} and Fig.~\ref{lock-in-perfect-pi-tri}, respectively.
To determine the lock-in frequency
$\omega_{\text{lock-in}} =
\omega_{\text{lock-in}}(\tau_2,K_{\rm vco}/\tau_1)$
by diagrams in Fig.~\ref{lock-in-perfect-pi-sin} and Fig.~\ref{lock-in-perfect-pi-tri}
one has to choose a curve corresponding to the loop filter parameter $\tau_2 \in [0,1]$,
select a point on the curve with the abscissa equal to $\frac{K_{\rm vco}}{\tau_1}$,
and then corresponding ordinate in the left subfigure
and product of the ordinates by $\frac{K_{\rm vco}}{\tau_1}$
in the right subfigure
gives the lock-in frequency $\omega_{\text{lock-in}}$
(see Fig.~\ref{lock-in-diagram-expl}).

Note, that for arbitrary $\tau_2>0$
by the following change of variables
\(
  t \to  t\tau_2, x \to x \tau_2
\)
in system \eqref{PLLPI-reduced-params}
we get
$\omega_{e}^{\rm{free}} \to \tau_2\omega_{e}^{\rm{free}},
K_{\rm{vco}} \to  \tau_2^2 K_{\rm{vco}}/\tau_1$
(the considered time reparametrisation
does not affect the cycle slipping property \eqref{eq-cs-sup}
of trajectories).
Therefore, $\omega_{\text{lock-in}}(\tau_2,K_{\rm vco})$ for $\tau_2>1$
can be computed from the diagrams in Fig.~\ref{lock-in-perfect-pi-sin} and Fig.~\ref{lock-in-perfect-pi-tri} by the formula:
\begin{equation}
  \omega_{\text{lock-in}}(\tau_2,K_{\rm vco}/\tau_1)=
  \omega_{\text{lock-in}}(1,K_{\rm{vco}}\tau_2^2/\tau_1)/\tau_2.
\end{equation}

In \cite{AleksandrovKLYY-2016-arXiv-sin,AleksandrovKLNYY-2016-IFAC}
a method for analytical computation of separatrices in model \eqref{PLLPI}
is developed, which allows us to get the following lock-in frequency estimate
\begin{equation}
\begin{aligned}
 & \omega_{\text{lock-in}}(\tau_2,K_{\rm vco}/\tau_1) \approx \\
 & \frac{1}{\sqrt{2}}
 \bigg(\frac{K_{\rm vco}}{\tau_1}\bigg)^{\frac{1}{2}}+
 \frac{1}{6}\frac{K_{\rm vco}}{\tau_1}+
 \frac{\tau_2^2(5 - 6\ln 2)}{\sqrt{2}\,36\tau_1}
 \bigg(\frac{K_{\rm vco}}{\tau_1}\bigg)^{\frac{3}{2}}.
\end{aligned}
\end{equation}

Fig.~\eqref{lock-in-pi-pp-fig} shows the phase portraits of model \eqref{PLLPI}
with phase detector characteristics \eqref{sinusoidal_pd}
when increasing the parameter $\omega_{\rm{ref}}$ in the computation of the lock-in range.

\begin{figure*}[h]
  \centering
  \includegraphics[width=0.95\linewidth]{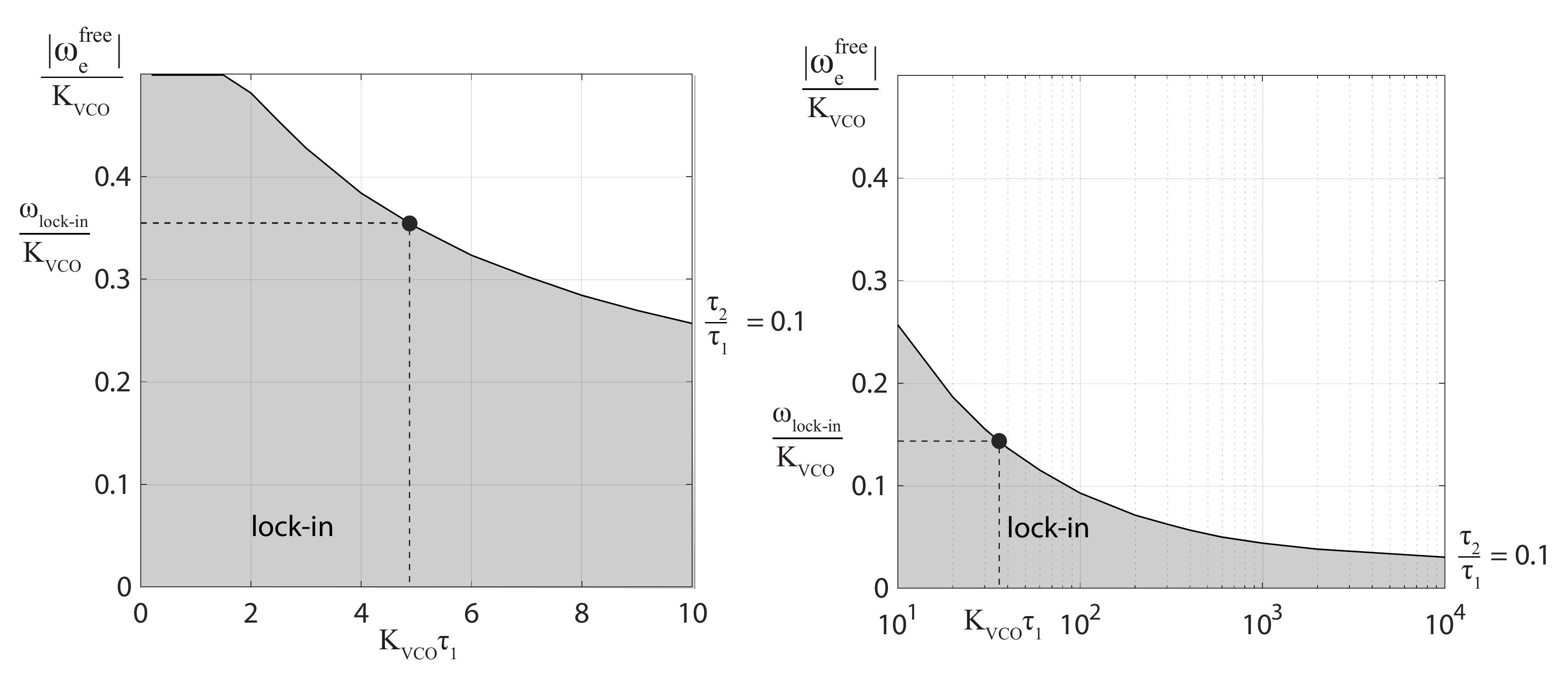}
  \caption{
  Lock-in range depending on $K_{\rm vco}\tau_1$ (shaded domains)
  for model \eqref{PLLleadlag} (lead-lag filter)
  with sinusoidal phase detector characteristics \eqref{sinusoidal_pd}
  and $\tau_2/\tau_1=0.1$.
  Left subfigure: $\frac{\omega_{\text{lock-in}}}{K_{\rm vco}}$
  for $0 < K_{\rm vco}\tau_1 < 10$
  in linear scale.
  Right subfigure: $\frac{\omega_{\text{lock-in}}}{K_{\rm vco}}$ for $K_{\rm vco}\tau_1 > 10$ in logarithmic scale.
  }
  \label{lock-in-lead-lag-expl}
\end{figure*}

\begin{figure*}[h]
  \centering
  \includegraphics[width=0.95\linewidth]{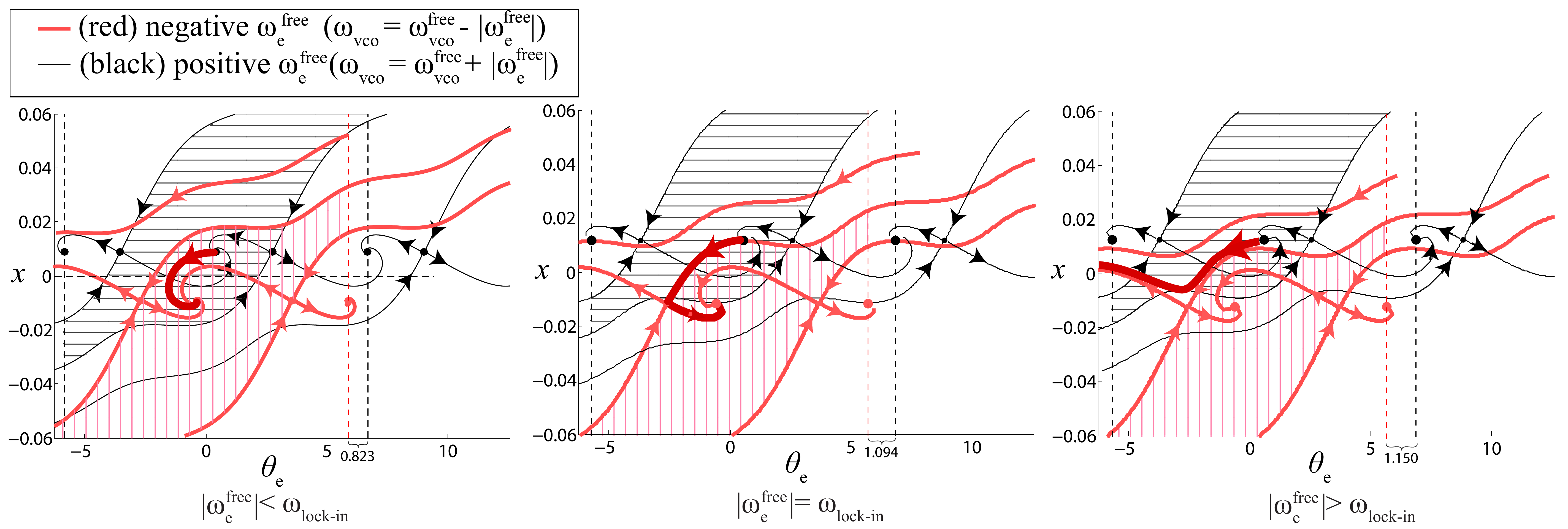}
  \caption{
  Phase portraits for model \eqref{PLLleadlag} with the following parameters:
 $H(s)= \frac{1+s\tau_2}{1+s\tau_1}$,
 $\tau_1 = 6.33\cdot10^{-2}$,
 $\tau_2 = 1.85\cdot10^{-2}$,
 $K_{\rm vco}=250$, $v_e(\theta_e)= \frac{1}{2}\sin(\theta_e)$.
 Black color is for the model with positive $\omega_e^{\text{free}}=|\widetilde{\omega}|$.
 Red is for the model with negative $\omega_e^{\text{free}}=-|\widetilde{\omega}|$.
 Equilibria (dots), separatrices pass in and out of the saddles equilibria,
 local lock-in domains are shaded
 (upper black horizontal lines are for $\omega_e^{\text{free}}>0$,
 lower red vertical lines are for $\omega_e^{\text{free}}<0$).
 Left subfig: $\omega_e^{\text{free}} = \pm 50$;
 middle subfig: $\omega_e^{\text{free}} = \pm 65$;
 right subfig: $\omega_e^{\text{free}} = \pm 68$.}
  \label{lock-in-lead-lag-pp-fig}
\end{figure*}

\subsection{Lead-lag filter}
Consider the PLL with passive lead-lag loop filter
having transfer function
$H(s) = \frac{1+s\tau_2}{1 + s\tau_1}$, where $\tau_1>\tau_2\geq 0$.
The model \eqref{final_system} becomes
\begin{equation}\label{PLLleadlag}
  \begin{aligned}
    & \dot x = -\frac{1}{\tau_1} x + \frac{v_e(\theta_e(t))}{\tau_1},\\
    & \dot\theta_{e} = \omega_{e}^{\rm{free}}
   - K_{\rm{vco}}
   \bigg(
    \Big(1-\frac{\tau_2}{\tau_1}\Big)x + \frac{\tau_2}{\tau_1}v_e(\theta_e(t))
   \bigg).
  \end{aligned}
\end{equation}
Model \eqref{PLLleadlag} with PD characteristic \eqref{sinusoidal_pd}
for $\omega_e^{\rm{free}}< K_{\rm vco}/2$ has one asymptotically stable
$\Big(\sin^{-1}\left(\frac{2\omega_e^{\rm{free}}}{K_{\rm vco}}\right), -\frac{\omega_e^{\rm{free}}}{K_{\rm vco}}\Big)$
and one unstable
$\Big(\pi~-~\sin^{-1}\left(\frac{2\omega_e^{\rm{free}}}{K_{\rm vco}}\right),-\frac{\omega_e^{\rm{free}}}{K_{\rm vco}}\Big)$;
and with PD characteristic \eqref{triangular_pd}
for $\omega_e^{\rm{free}}< K_{\rm vco}$
has one asymptotically stable equilibrium
$\Big(\frac{\pi}{2}\frac{\omega_e^{\rm{free}}}{K_{\rm vco}},-\frac{\omega_e^{\rm{free}}}{K_{\rm vco}}\Big)$
and one unstable equilibrium
$\Big(\frac{\pi}{2}\left( \frac{\omega_e^{\rm{free}}}{K_{\rm vco}}-2 \right),-\frac{\omega_e^{\rm{free}}}{K_{\rm vco}}\Big)$.
Thus, the pull-in range of model \eqref{PLLleadlag}
loop filter is bounded.
It can be estimated by the phase plane analysis or numerical computation\footnote{
By the Lyapunov function $V(x) = \frac{\tau_1}{2}x^2$
we get $\dot V(x) = -x^2 + v_{e}(\theta_e)x < 0$
for $|x| > v_{e}(\theta_e)$,
thus,
$\limsup_{t\to\infty}|x(t)|\leq \max v_{e}(\theta_e)=1$
and only a bounded domain with respect to $x$ has to be studied.
}
(see corresponding diagram in \cite{Belyustina-1970-eng,LeonovK-2014-book}),
one of the difficulties of its computation caused
by the so-called hidden oscillations \cite{LeonovK-2013-IJBC,LeonovKYY-2015-TCAS,KuznetsovLYY-2017-CNSNS}.
The computed 
\begin{figure}[H]
 \centering
 \includegraphics[width=0.85\linewidth]{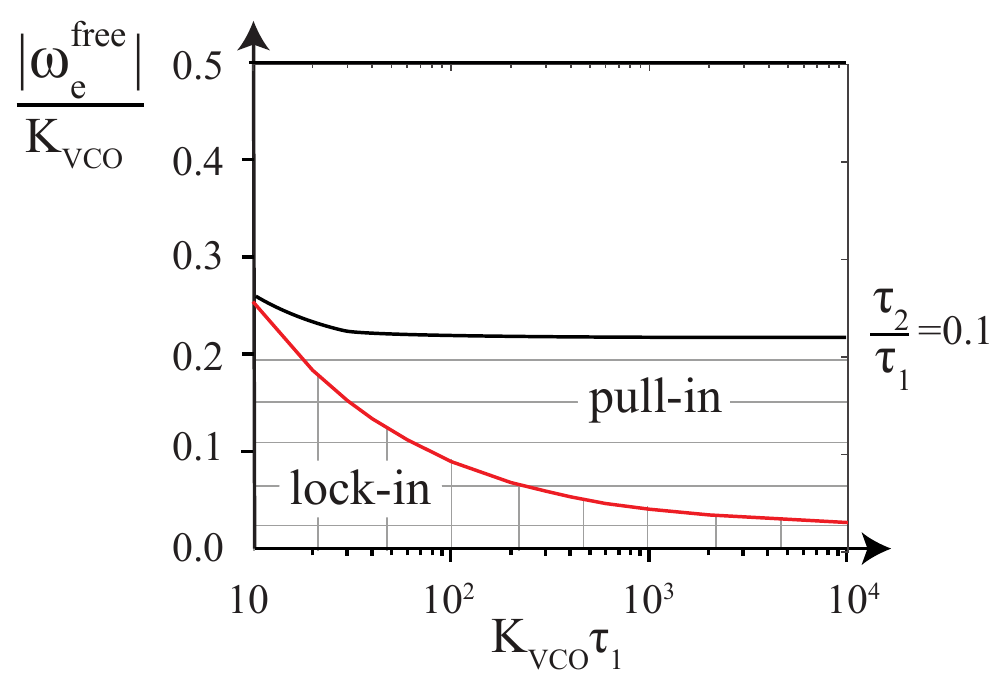}
 \caption{
  Pull-in and lock-in ranges depending on $K_{\rm vco}\tau_1$ (shaded domains)
  for model \eqref{PLLleadlag} (lead-lag filter)
  with sinusoidal phase detector characteristics \eqref{sinusoidal_pd}
  and $\tau_2/\tau_1=0.1$.
  }
 \label{ranges-comparison}
\end{figure}
\noindent above values of the lock-in frequency are less than corresponding
values of the pull-in frequency (see Fig.~\ref{ranges-comparison}).

 For an arbitrary $\tau_1>0$ the time reparametrization
 \(
    t \to t\tau_1>0
 \)
 reduces system \eqref{PLLleadlag} to the form
\begin{equation}\label{PLLleadlag-reduced-params}
  \begin{aligned}
    & \dot x = -x + v_e(\theta_e),\\
    & \dot\theta_{e} = \tau_1\omega_{e}^{\rm{free}}
   - \tau_1K_{\rm{vco}}
   \bigg(
    \Big(1-\frac{\tau_2}{\tau_1}\Big)x + \frac{\tau_2}{\tau_1}v_e(\theta_e(t))
   \bigg),
  \end{aligned}
\end{equation}
where $\tilde \omega_{e}^{\rm{free}} = \tau_1\omega_{e}^{\rm{free}}$ depends
on two parameters: $\tau_1K_{\rm{vco}}$ and $\tau_2/\tau_1$.

For model \eqref{PLLleadlag}
with phase detector characteristics \eqref{sinusoidal_pd} and \eqref{triangular_pd}
the lock-in frequency diagrams
as function of $\tau_2/\tau_1$ (for $\tau_2/\tau_1 \in [0,1]$) and $K_{\rm vco}\tau_1$
are shown in Fig.~\ref{lock-in-diagram-lead-lag-sin} and Fig.~\ref{lock-in-diagram-lead-lag-tri}, respectively.
To determine the lock-in frequency
$\omega_{\text{lock-in}} =
\omega_{\text{lock-in}}(\tau_2/\tau_1,K_{\rm vco}\tau_1)$
by diagrams in Fig.~\ref{lock-in-diagram-lead-lag-sin} and Fig.~\ref{lock-in-diagram-lead-lag-tri}
one has to choose a curve corresponding to the loop filter parameters $\tau_2/\tau_1 \in [0,1]$,
select a point on the
curve
with abscissa equal to $K_{\rm vco}\tau_1$,
and then the product of the ordinates by $K_{\rm vco}$
gives the lock-in frequency $\omega_{\text{lock-in}}$
(see Fig.~\ref{lock-in-lead-lag-expl}).

Fig.~\eqref{lock-in-lead-lag-pp-fig} shows phase portraits of model \eqref{PLLleadlag}
with phase detector characteristics \eqref{sinusoidal_pd}
when increasing the parameter $\omega_{\rm{ref}}$ in the computation of the lock-in range.


\section*{Conclusion}Various PLL-based circuits
(see, e.g.
optical Costas loop used in intersatellite communication,
BPSK Costas loop, two-phase Costas loop, two-phase PLL and others
\cite{emura2000high,BestKLYY-2014-IJAC,Best-2007,LeonovKYY-2015-SIGPRO,BestKKLYY-2015-ACC,rosenkranz2016receiver,BestKLYY-2016,KuznetsovLYY-2017-CNSNS,Middlestead-2017})
are represented in the signal's phase space by
the model in Fig.~\ref{phase-space-model-fig}
and, thus,
their lock-in ranges can be estimated by the above diagrams.
Remark that the change of variables
$(\theta_e,x) \to (\theta_e/K_{\rm p}, K_{\rm pd}x)$
transforms the model \eqref{final_system}
with PD characteristic $K_{\rm pd}v_e\big(K_{p}\theta_e(t)\big)$
to the from with PD characteristics $v_e(\theta_e(t))$
(i.e. with $K_{\rm p} \to 1, K_{\rm pd} \to 1$)
and
$\omega_e^{\rm free} \to K_{\rm p}\omega_e^{\rm free},
K_{\rm vco} \to K_{\rm p}K_{\rm pd}K_{\rm vco}$.
Thus, the lock-in frequency for the model with active PI loop filter
can be computed from the diagrams in Fig.~\ref{lock-in-perfect-pi-sin} and Fig.~\ref{lock-in-perfect-pi-tri}
as $\omega_{\text{lock-in}}(\tau_2,K_{\rm p}K_{\rm pd}K_{\rm vco}/\tau_1)/K_{\rm p}$
and for the model with lead-lag loop filter
can be computed from diagrams in Fig.~\ref{lock-in-diagram-lead-lag-sin} and Fig.~\ref{lock-in-diagram-lead-lag-tri}
as
$\omega_{\text{lock-in}}(\tau_2/\tau_1,K_{\rm p}K_{\rm pd}K_{\rm vco}\tau_1)/K_{\rm p}$.
Analytical estimates
of the lock-in range for two-dimensional models
can be obtained by the Andronov point-transformation method \cite{AndronovVKh-1937}
and the study of separatrices in cylindrical phase space
(some estimations of separatrices for the classical PLL
can be found in
\cite{Gubar-1961,Shakhtarin-1969,ShahgildyanL-1966,Shalfeev-2013-book,Gardner-2005-book,HuqueS-2013}).
For the second-order PLL with PI filter and triangular phase-detector characteristic
an analytical estimate of the lock-in range can be found in
\cite{AleksandrovKLNYY-2016-IFAC,AleksandrovKLYY-2016-arXiv-sin,AleksandrovKLYY-2016-arXiv-impulse}.


\section*{Acknowledgments}
This work was supported by the Russian Science Foundation (14-21-00041).
The authors would like to thank Roland Best,
the founder of the Best Engineering Company (Oberwil, Switzerland)
and the author of the bestseller on PLL-based circuits \cite{Best-2007}
for valuable discussion on the lock-in range concept.

\bibliographystyle{IEEEtran}


\newpage
\begin{IEEEbiography}
 {Nikolay Kuznetsov}
 received his Candidate degree from Saint-Petersburg State University
 (2004), Ph.D. from the University of Jyv\"{a}skyl\"{a} (2008),
 and D. Sc. from Saint-Petersburg State University (2016).
 He is currently Professor and Deputy Head of the Department of Applied Cybernetics
 at Saint-Petersburg State University,
 Visiting Professor at the University of Jyv\"{a}skyl\"{a},
 and member of the MERLIN international research group at the Ton Duc Thang University.
 His interests are now in dynamical systems stability and oscillations,
 Lyapunov exponent, chaos, hidden attractors, phase-locked loop nonlinear analysis,
 nonlinear control systems.
 E-mail: nkuznetsov239@gmail.com (corresponding author)
\end{IEEEbiography}

\begin{IEEEbiography}
{Gennady Leonov}
received his Candidate degree in 1971 and D. Sc. in 1983 from Saint-Petersburg State University.
In 1986 he was awarded the USSR State Prize
\emph{for development of the theory of phase synchronization for radiotechnics and communications}.
Since  1988 he has been Dean of the Faculty of Mathematics and Mechanics
at Saint-Petersburg State University
and since 2007 Head of the Department of Applied Cybernetics.
He is member (corresponding) of the Russian Academy of Science, in 2011 he was elected to the IFAC Council.
His research interests are now in control theory and dynamical systems.
\end{IEEEbiography}

\begin{IEEEbiography}
{Marat Yuldashev}
received his Candidate degree from St.Petersburg State University
(2013) and Ph.D. from the University of Jyv\"{a}skyl\"{a} (2013).
He is currently at Saint-Petersburg University.
His research interests cover nonlinear models of phase-locked loops and Costas loops, and SPICE simulation.
\end{IEEEbiography}

\begin{IEEEbiography}
{Renat Yuldashev}
received his Candidate degree from St.Petersburg State University
(2013) and Ph.D. from the University of Jyv\"{a}skyl\"{a} (2013).
He is currently at Saint-Petersburg University.
His research interests cover nonlinear models of phase-locked loops and Costas loops,
and simulation in MatLab Simulink.
\end{IEEEbiography}

\end{document}